\documentclass[12pt,a4paper]{article}
\usepackage{german}
\usepackage{amssymb}
\usepackage{amsmath}
\usepackage{graphicx}
\usepackage{latexsym}
\usepackage{vmargin}
\usepackage[latin1]{inputenc}

\newcommand{\R}{\mbox{I}\hspace{-0.05cm}\mbox{R}}

\newcommand{\C}{\mathbb{C}}
\begin{document}
\newtheorem{theorem}{Satz}[section]
\newtheorem{remark}{Bemerkung}[section]
\newtheorem{aufgabe}{Aufgabe}[section]
\newtheorem{define}{Definition}[section]
\newtheorem{example}{Beispiel}[section]
\newtheorem{lemma}{Lemma}[section]
\begin{center}
{\Large \bf \"Uber die Anwendung des Tschebyschew-Verfahrens zum Ausbau des Weierstra\ss-Kerner-Verfahrens}\\[2ex]
{\bf Uwe Sch\"afer}\\[1ex]
Institut f\"ur Angewandte und Numerische Mathematik,\\
Karlsruher Institut f\"ur Technologie\\
Kaiserstra\ss e 12, D-76128 Karlsruhe, Germany\\
E-Mail: Uwe.Schaefer@kit.edu\\[4ex]
{\it F\"ur Prof. Dr. G. Alefeld, der mich zu dieser Arbeit ermunterte.}
\end{center}
\begin{center}
{\bf Abstract}
\end{center}
We extend the Weierstra\ss-Kerner method by applying the Chebychev method to the function $F$ that Kerner has used to show that the formula of Weierstra\ss { }actually is the Newton method applied to that $F$. The resulting method is already known but we want to present the process in one go and in a detailed way.
\section{Einleitung}
\label{intro}
Das Weierstra\ss-Kerner-Verfahren ist ein Verfahren zur simultanen Berechnung der Nullstellen von Polynomen. W\"ahrend Weierstra\ss { }dieses Verfahren bereits 1891 benutzte, um einen damals neuen Beweis des Fundamentalsatzes der Algebra zu f\"uhren, bewies Kerner 1966, dass das Verfahren in Wirklichkeit das Newton-Verfahren ist angewandt auf eine Funktion $F:\C^n \to \C^n$, die aus den elementarsymmetrischen Polynomen aus der Satzgruppe von Vieta besteht. 

Dieser Artikel entwickelt Kerners Ansatz weiter, indem auf $F$ das Tschebyschew-Verfahren angewandt wird. Die resultierenden Formeln wurden allerdings bereits 1983  (auf einem anderen Weg) von Tanabe gefunden, und in \cite{kanno} wurde auch bereits gezeigt, dass das Verfahren von Tanabe genau das Tschebyschew-Verfahren angewandt auf $F$ ist. Wir wollen in diesem Artikel die Herleitung dennoch einmal detailliert durchexerzieren, damit die ganze Entstehung der Formeln in einem Guss pr\"asentiert wird.  
\subsection{Polynome}
Wir betrachten normierte reelle  Polynome vom Grad $n$, d.h. $p:\R \to \R$ mit 
\begin{equation} \label{polynom}
p(t)= t^n + a_{n-1} t^{n-1} + \cdots + a_{1} t + a_{0}, \quad a_{j} \in \R, \, j = 0,...,n-1.
\end{equation}
Der Fundamentalsatz der Algebra besagt, dass jedes Polynom $n$-ten Grades $n$ (nicht notwendigerweise paarweise verschiedene und wom\"oglich komplexe) Nullstellen besitzt. Es existieren also
\begin{equation} \label{nullen}
\xi_{j}\in \C,\quad  \, j=1,...,n \quad \mbox{mit} \quad p(\xi_j)=0,\quad  \, j=1,...,n,
\end{equation}
und $p$ l\"asst sich darstellen durch
\[
p(t)=  \prod \limits_{j=1}^{n} (t-\xi_{j}).
\]
Bekanntlich gibt es nur f\"ur $n\leq4$ explizite Formeln f\"ur die Nullstellen, siehe \cite{artin}, \cite{bew}, \cite{froeba}. Daher muss man sich mit N\"aherungsverfahren begn\"ugen. Allerdings gibt es durch Koeffizientenvergleich einen einfachen Zusammenhang zwischen den Nullstellen $\xi_j$ und den Koeffizienten $a_j$.
\begin{theorem} (Satz von Vieta)
Es gilt
\[
a_{n-k}= \sum \limits_{i_{1} < i_{2} < \cdots <i_{k}}
(- \xi_{i_{1}}) \cdot (- \xi_{i_{2}} ) \cdots  (- \xi_{i_{k}})\, , \quad k=n,...,1.
\]
\end{theorem}
Beispiel $n=4$: Dann ist
\[
p(t)=(t-\xi_{1}) \cdot(t-\xi_{2}) \cdot  (t-\xi_{3}) \cdot (t-\xi_{4}).
\]
Der Satz von Vieta erspart das Ausmultiplizieren und man erh\"alt sofort
\[
\begin{array}{rcl}
p(t) & = & t^4 + t^3(-\xi_{1}-\xi_{2}-\xi_{3}-\xi_{4}) \\[1ex]
     & + & t^2 (\xi_{1} \xi_{2} + \xi_{1} \xi_{3} + \xi_{1} \xi_{4} + \xi_{2}\xi_{3}
            +\xi_{2} \xi_{4} +\xi_{3} \xi_{4})  \\[1ex]
& + & t(-\xi_{1}\xi_{2} \xi_{3} -\xi_{1}\xi_{2}\xi_{4} - \xi_{1}\xi_{3}\xi_{4}
-\xi_{2}\xi_{3}\xi_{4} ) \\[1ex]
& + & \xi_{1}\xi_{2} \xi_{3}\xi_{4} .
\end{array}
\]
Kerners erste Idee war es, eine Funktion $V =(v_{i}): \C^n \to \C^n$ durch
\[
v_{n-\nu +1}({\bf x})=\sum \limits_{i_{1} < i_{2} < \cdots <i_{\nu}}
(- x_{i_{1}}) \cdot (- x_{i_{2}} ) \cdots  (- x_{i_{\nu}}) , \quad \nu=n,...,1
\]
zu definieren.\\[1ex]
Beispiel $n=4$: Dann ist
\[
\begin{array}{rcl}
v_{1}({\bf x}) & = & x_{1}  x_{2} x_{3} x_{4}\\[1ex]
v_{2}({\bf x}) & = & -x_{1}x_{2} x_{3} -x_{1}x_{2}x_{4} - x_{1}x_{3}x_{4}
-x_{2}x_{3}x_{4} \\[1ex]
v_{3}({\bf x}) & = & x_{1} x_{2} + x_{1} x_{3} + x_{1} x_{4} + x_{2}x_{3}
            +x_{2} x_{4} +x_{3} x_{4} \\[1ex]
v_{4}({\bf x}) & = & -x_{1}-x_{2}-x_{3}-x_{4} .
\end{array}
\]
Ist das Polynom (\ref{polynom}) gegeben durch den
Vektor
\[
{\bf a}=\left(
\begin{array}{c}
a_{0} \\ a_{1} \\ \vdots \\ a_{n-1}
\end{array}
\right),
\]
so ist der Vektor der Nullstellen aus (\ref{nullen})
\[
{\bf \xi}=\left(
\begin{array}{c}
\xi_{1} \\ \xi_{2} \\ \vdots \\ \xi_{n}
\end{array}
\right)
\]
die Nullstelle der Funktion $F:\C^n \to \C^n$ mit
\begin{equation} \label{funktionF}
F({\bf x}):= V({\bf x}) - {\bf a}.
\end{equation}
Kerners zweite Idee war es, auf diese Funktion das Newton-Verfahren anzuwenden: 
\begin{equation} \label{newton}
\left.
\begin{array}{rcl}
{\bf x^{(0)}} &  \in \C^n &\mbox{ beliebig} \\[1ex]
 {\bf x^{(m+1)}} & := & {\bf x^{(m)}} - \Big( F'({\bf x^{(m)}})\Big)^{-1} \cdot F({\bf x^{(m)}}).
\end{array}
\right\}
\end{equation}
{\bf Beobachtung 1} F\"ur festes $ {\bf x^{(m)}} \in \C^n$ setzen wir
\begin{equation} \label{vektorb}
 {\bf b} = \left(
\begin{array}{c}
b_{0} \\ b_{1} \\ \vdots \\ b_{n-1}
\end{array}
\right) := V({\bf x^{(m)}}).
\end{equation}
Nach dem Satz von Vieta gilt dann mit $b_{n}:=1$
\begin{equation} \label{beob1}
\prod \limits_{j=1}^{n} (t-x^{(m)}_{j}) = \sum \limits_{j=0}^{n} b_{j} t^j.
\end{equation}
{\bf Beobachtung 2} F\"ur $k=1,...,n$ gilt
\begin{equation} \label{genial}
\Big(1 \,\,  t \,\,  t^2 \,\,  \cdots \,\,  t^{n-1} \Big) \cdot F'({\bf x^{(m)}})_{\cdot k}
=  -\prod \limits_{j=1,j\not=k}^{n} (t-x^{(m)}_{j}).
\end{equation}
{\it Beweis:} Es ist
\begin{equation} \label{wieF'}
F'({\bf x})_{\cdot k} =
\left(
\begin{array}{c}
\displaystyle{\frac{\partial v_{1}({\bf x})}{\partial x_{k}}} \\
\vdots \\
\displaystyle{\frac{\partial v_{n}({\bf x})}{\partial x_{k}}}
\end{array}
\right) =
\left(
\begin{array}{c}
\displaystyle{-\sum \limits_{i_{1} < i_{2} < \cdots <i_{n-1}} \, \,
\prod \limits_{j=1,\, i_{j} \not=k}^{n-1}  (- x_{i_{j}} ) }  \\[3ex]
\displaystyle{-\sum \limits_{i_{1} < i_{2} < \cdots <i_{n-2}} \, \,
\prod \limits_{j=1,\, i_{j} \not=k}^{n-2}  (- x_{i_{j}} ) } \\[2ex]
\vdots \\
\displaystyle{-\sum \limits_{i=1,\, i\not=k}^{n}
(- x_{i} ) } \\[3ex]
-1
\end{array}
\right).
\end{equation}
Beispiel $n=4$: Es ist 
\[
F'({\bf x})_{\cdot 1}= \left(
\begin{array}{c}
x_{2} x_{3} x_{4} \\
-x_{2} x_{3} - x_{2} x_{4} -x_{3} x_{4} \\
  x_{2} + x_{3} + x_{4} \\
  -1
\end{array}
\right), \,
F'({\bf x})_{\cdot 2}= \left(
\begin{array}{c}
x_{1} x_{3} x_{4} \\
-x_{1} x_{3} - x_{1} x_{4} -x_{3} x_{4} \\
  x_{1} + x_{3} + x_{4} \\
  -1
\end{array}
\right)
\]
und
\[
F'({\bf x})_{\cdot 3}= \left(
\begin{array}{c}
x_{1} x_{2} x_{4} \\
-x_{1} x_{2} - x_{1} x_{4} -x_{2} x_{4} \\
  x_{1} + x_{2} + x_{4} \\
  -1
\end{array}
\right), \,
F'({\bf x})_{\cdot 4}= \left(
\begin{array}{c}
x_{1} x_{2} x_{3} \\
-x_{1} x_{2} - x_{1} x_{3} -x_{2} x_{3} \\
  x_{1} + x_{2} + x_{3} \\
  -1
\end{array}
\right).
\]
Mit ${\bf x}={\bf x^{(m)}}$ erh\"alt man mit der Formel aus dem Satz von Vieta, dass in der $k$-ten
Spalte von $F'({\bf x^{(m)}}) $ die mit $-1$
multiplizierten Koeffizienten des normierten Polynoms vom Grad $n-1$ mit den Nullstellen
\[
x_{1}^{(m)},\cdots , \, x_{k-1}^{(m)} , x_{k+1}^{(m)}, \cdots , \, x_{n}^{(m)}
\]
stehen. Dies \"aquivalent umgeschrieben ist (\ref{genial}).\hfill $\Box$\\[2ex]
Mit Hilfe von (\ref{genial}) ist es sehr einfach, die Inverse von
$F'({\bf x^{(m)}}) $ explizit anzugeben.\\[2ex]
{\bf Beobachtung 3} Es gelte
\[
 x_{i}^{(m)} \not= x_{j}^{(m)}, \quad i\not=j.
\]
Dann ist $F'({\bf x^{(m)}}) $ invertierbar, und es gilt mit
\begin{equation} \label{wfunk}
B^{(k)}(t):= -\prod \limits_{j=1,j\not=k}^{n} (t-x^{(m)}_{j}), \quad k=1,...,n
\end{equation}
folgende Gleichheit
\[
 F'({\bf x^{(m)}})^{-1}
 =\left(
 \begin{array}{c}
 \displaystyle{\frac{1}{B^{(1)}(x_{1}^{(m)})}}
 \left(
 1 \quad x_{1}^{(m)} \quad  (x_{1}^{(m)})^2 \cdots  (x_{1}^{(m)})^{n-1}
 \right)\\[3ex]
 \vdots\\[3ex]
\displaystyle{\frac{1}{B^{(n)}(x_{n}^{(m)})}}
 \left(
 1 \quad x_{n}^{(m)} \quad  (x_{n}^{(m)})^2 \cdots  (x_{n}^{(m)})^{n-1}
 \right)
 \end{array}
 \right).
\]
{\it Beweis:} Wegen (\ref{genial}) gilt
\[
 \displaystyle{\frac{1}{B^{(l)}(x_{l}^{(m)})}}
 \left(
 1 \quad x_{l}^{(m)} \quad  (x_{l}^{(m)})^2 \cdots  (x_{l}^{(m)})^{n-1}
 \right) \cdot F'({\bf x^{(m)}})_{\cdot k} =
 \frac{B^{(k)}(x_{l}^{(m)})}{B^{(l)}(x_{l}^{(m)})} = \delta_{lk}.
\]
$\delta_{lk}$ bezeichnet dabei das Kronecker-Symbol. \hfill $\Box$\\[2ex]
Somit ergibt sich f\"ur das Newton-Verfahren \"uber (\ref{vektorb}) und
(\ref{beob1}) und da $a_{n}=1$, $b_{n}=1$ f\"ur $l=1,...,n$
\[
\begin{array}{rcl}
x_{l}^{(m+1)} & = &
x_{l}^{(m)} - \displaystyle{\frac{1}{B^{(l)}(x_{l}^{(m)})}}
 \left(
 1 \quad x_{l}^{(m)} \quad  (x_{l}^{(m)})^2 \cdots  (x_{l}^{(m)})^{n-1}
 \right) \left(
 \begin{array}{c}
 b_{0} - a_{0} \\ \vdots \\ b_{n-1} - a_{n-1}
  \end{array}
 \right)\\[3ex]
 & = &  x_{l}^{(m)} - \displaystyle{\frac{1}{B^{(l)}(x_{l}^{(m)})}
 \sum \limits_{j=0}^{n}} (b_{j}-a_{j}) \cdot  (x_{l}^{(m)})^j \\[3ex]
 & = & x_{l}^{(m)} + \displaystyle{\frac{1}{B^{(l)}(x_{l}^{(m)})}
 \sum \limits_{j=0}^{n}} a_{j} \cdot (x_{l}^{(m)})^j  =
 x_{l}^{m} + \displaystyle{\frac{p(x_{l}^{(m)})}{B^{(l)}(x_{l}^{(m)})}} \\[3ex]
 & = & x_{l}^{(m)} - \displaystyle{\frac{p(x_{l}^{(m)})}{\prod \limits_{j=1,j\not=l}^{n}
 (x_{l}^{(m)}-x_{j}^{(m)})}} .
 \end{array}
\]
Die letzte Formel ist die Iterationsvorschrift, die man gerne Weierstra\ss-Kerner-Verfahren nennt. Der Vollst\"andigkeit wegen wollen wir betonen, dass es auch andere N\"aherungsverfahren gibt. Wir verweisen auf \cite{alefeld}, \cite{kyurk} und \cite{petkovic}.\\[2ex]
Bevor wir im Hauptteil das Tschebyschew-Verfahren auf $F$ aus (\ref{funktionF}) anwenden, wollen wir noch in einem Abschnitt die zweite Ableitung als Bilinearform wiederholen (siehe \cite{wloka}), da wir sie beim Tschebyschew-Verfahren benutzen.
\subsection{Die zweite Ableitung als Bilinearform}
Es sei $F=(f_{i}) : \C^n \to \C^n$, d.h. $f_{i} :
\C^n \to \C$ f\"ur $i=1,...,n$. Es sei
$F$ zweimal stetig differenzierbar. Dann gilt
\[
F'({\bf x}) =
\left(
\begin{array}{c}
\mbox{grad }f_{1}({\bf x}) \\
\vdots \\
\mbox{grad }f_{n}({\bf x})
\end{array}
\right)
\]
und
\begin{equation} \label{bilinear}
F''({\bf x}) =
\left(
\begin{array}{c|c|c|c}
\mbox{grad }\displaystyle{\frac{ \partial f_{1}({\bf x})}{\partial x_{1}}} &
\mbox{grad }\displaystyle{\frac{ \partial f_{1}({\bf x})}{\partial x_{2}}} &
\cdots \cdots&\mbox{grad }\displaystyle{\frac{ \partial f_{1}({\bf x})}{\partial x_{n}}} \\[2ex]
\vdots &
\cdots & \cdots \cdots & \vdots \\[2ex]
\mbox{grad }\displaystyle{\frac{ \partial f_{n}({\bf x})}{\partial x_{1}}}&
\mbox{grad }\displaystyle{\frac{ \partial f_{n}({\bf x})}{\partial x_{2}}} &\cdots  \cdots&
\mbox{grad }\displaystyle{\frac{ \partial f_{n}({\bf x})}{\partial x_{n}}}
\end{array}
\right).
\end{equation}
In der $i$-ten Zeile von $F''({\bf x})$ steht somit zeilenweise hintereinander
geschrieben die Hesse-Matrix von $f_{i}({\bf x})$. Man kann auch schreiben
\[
 F''({\bf x})= \left(
\begin{array}{c|c|c|c}
A^{(1)}({\bf x}) & A^{(2)}({\bf x}) &
\cdots & A^{(n)}({\bf x})
\end{array}
\right)
\]
mit
\[
A^{(k)}({\bf x}) =(a^{(k)}_{ij}({\bf x}))\in \C^{n \times n}, k=1,...,n.
\]
$F''({\bf x})$ ist eine so genannte Bilinearform. Es gilt f\"ur ${\bf y}, {\bf z}\in \C^n$
\[
F''({\bf x})({\bf y}, {\bf z})= \underbrace{\left(
\begin{array}{c|c|c|c}
A^{(1)}({\bf x})\cdot {\bf y} & A^{(2)}({\bf x})\cdot {\bf y} &
\cdots & A^{(n)}({\bf x})\cdot{\bf y}
\end{array}
\right)}_{\in \C^{n \times n}}\cdot  {\bf z} \in \C^n.
\]
Es gilt
\[
F''({\bf x})({\bf y}, {\bf z})_{i}= \sum \limits_{k=1}^{n}
\Big( \sum \limits_{j=1}^{n} a^{(k)}_{ij}({\bf x})\cdot y_{j}\Big) z_{k} , \quad i=1,...,n.
\]
Aufgrund des Satzes von Schwarz gilt $a^{(k)}_{ij}({\bf x})=a^{(j)}_{ik}({\bf x})$ und somit
\[
F''({\bf x})({\bf y}, {\bf z})= F''({\bf x})({\bf z}, {\bf y}).
\]
Dies wird im Folgenden aber keine Rolle spielen, da in unserer Anwendung ${\bf y}= {\bf z} $ gelten wird.
\begin{lemma} Es seien $B \in \C^{n \times n}$, ${\bf y},\, {\bf z} \in \C^n$. Dann gilt
\[
B \cdot \Big( F''({\bf x})({\bf y}, {\bf z}) \Big) =\left(
\begin{array}{c|c|c|c}
B \cdot A^{(1)}({\bf x}) & B \cdot A^{(2)}({\bf x}) &
\cdots & B \cdot A^{(n)}({\bf x})
\end{array}
\right)({\bf y}, {\bf z}).
\]
\end{lemma}
{\it Beweis:} Es ist
\[
\begin{array}{rcl}
B \cdot \Big( F''({\bf x})({\bf y}, {\bf z}) \Big)&=&
B \cdot
\left(
\begin{array}{c|c|c|c}
A^{(1)}({\bf x})\cdot{\bf y} & A^{(2)}({\bf x})\cdot{\bf y} &
\cdots & A^{(n)}({\bf x})\cdot{\bf y}
\end{array}
\right)\cdot{\bf z}\\[2ex]
& = & \left(
\begin{array}{c|c|c|c}
B\cdot A^{(1)}({\bf x})\cdot {\bf y} & B \cdot A^{(2)}({\bf x})\cdot {\bf y} &
\cdots & B\cdot A^{(n)}({\bf x}) \cdot {\bf y}
\end{array}
\right)\cdot {\bf z}
\end{array}
\]
wegen der Assoziativit\"at der
Matrizen/Vektor-Multiplikation. \hfill $\Box$
\section{Das Tschebyschew-Verfahren angewandt auf $F$}
Es sei $F=(f_{i}) : \C^n \to \C^n$, d.h. $f_{i} :
\C^n \to \C$ f\"ur $i=1,...,n$, und es sei
$F$ zweimal stetig differenzierbar. Falls $F({\bf x^*})=o$, so erh\"alt man aus der
Taylorformel analog zur Herleitung
des Newton-Verfahrens
\[
 o = F({\bf x^*}) \approx F({\bf x^{(m)}}) + F'({\bf x^{(m)}})\cdot({\bf x^*}- {\bf x^{(m)}})
 + \frac{1}{2!} F''({\bf x^{(m)}})({\bf x^*}- {\bf x^{(m)}},{\bf x^*}- {\bf x^{(m)}}) .
\]
Man erh\"alt \"uber
\[
 {\bf x^*}- {\bf x^{(m)}} \approx - F'({\bf x^{(m)}})^{-1} \cdot  F({\bf x^{(m)}})
\]
das so genannte Tschebyschew-Verfahren: $ {\bf x^{(0)}} \in \C^n$ beliebig und
\[
{\bf x^{(m+1)}}=
 {\bf x^{(m)}} -
 \]
 \[
 F'({\bf x^{(m)}})^{-1} \cdot
 \left(
 F({\bf x^{(m)}}) +  \frac{1}{2} F''({\bf x^{(m)}})\Big(
 F'({\bf x^{(m)}})^{-1}   F({\bf x^{(m)}}),F'({\bf x^{(m)}})^{-1}   F({\bf x^{(m)}})\Big)
 \right).
\]
Es ist
\[
F''({\bf x})= \left(
\begin{array}{c|c|c|c}
A^{(1)}({\bf x}) & A^{(2)}({\bf x}) &
\cdots & A^{(n)}({\bf x})
\end{array}
\right)
\]
mit
\[
A^{(k)}({\bf x})=
\left(
\begin{array}{c}
\mbox{grad} \displaystyle{\frac{\partial f_{1}({\bf x})}{\partial x_{k}}} \\[1ex]
\vdots \\[1ex]
 \mbox{grad} \displaystyle{\frac{\partial f_{n}({\bf x})}{\partial x_{k}}}
\end{array}
\right).
\]
F\"ur (\ref{funktionF}) ergibt sich mit (\ref{wieF'}) f\"ur $l \not=k$
\[
\Big( A^{(k)}({\bf x}) \Big)_{\cdot l} =
\left(
\begin{array}{c}
\displaystyle{\frac{\partial^2 f_{1}({\bf x})}{\partial x_{k} \partial x_{l}}} \\[1ex]
\vdots \\[1ex]
\displaystyle{\frac{\partial^2 f_{n}({\bf x})}{\partial x_{k} \partial x_{l}}}
\end{array}
\right) =
\left(
\begin{array}{c}
\displaystyle{\sum \limits_{i_{1}  < \cdots <i_{n-2}} \, \,
\prod \limits_{\begin{array}{c} j=1,\\ i_{j} \not=k, \\ i_{j}\not=l
\end{array}}^{n-2}  (- x_{i_{j}} ) }
\\[10ex]
\displaystyle{\sum \limits_{i_{1}  < \cdots <i_{n-3}} \, \,
\prod \limits_{\begin{array}{c} j=1,\\ i_{j} \not=k, \\ i_{j}\not=l
\end{array}}^{n-3}  (- x_{i_{j}} ) }
\\[2ex]
\vdots \\
\displaystyle{\sum \limits_{i=1,i\not= k, i \not= l}^{n} (- x_{i})}\\
 1 \\
 0
\end{array}
\right)
\]
und
\[
 \Big( A^{(k)}({\bf x}) \Big)_{\cdot l}=o \quad \mbox{f\"ur }l=k.
\]
Beispiel $n=4$: Hier ergibt sich
\[
A^{(1)}({\bf x})=
\left(
\begin {array}{cccc}
0 & x_{3}x_{4} &   x_{2}x_{4} &  x_{2}x_{3} \\[1ex]
0 & -x_{3}-x_{4} &   -x_{2}-x_{4} &  -x_{2}-x_{3} \\[1ex]
0 & 1 & 1 & 1 \\[1ex]
0 & 0 & 0 & 0
\end{array}
\right)
\]
\[
A^{(2)}({\bf x})=
\left(
\begin {array}{cccc}
x_{3}x_{4} & 0 &   x_{1}x_{4} &  x_{1}x_{3} \\[1ex]
-x_{3}-x_{4} &  0 &  -x_{1}-x_{4} &  -x_{1}-x_{3} \\[1ex]
1 & 0 & 1 & 1 \\[1ex]
0 & 0 & 0 & 0
\end{array}
\right)
\]
 \[
A^{(3)}({\bf x})=
\left(
\begin {array}{cccc}
x_{2}x_{4} &   x_{1}x_{4} & 0 & x_{1}x_{2} \\[1ex]
-x_{2}-x_{4} &   -x_{1}-x_{4} & 0 &  -x_{1}-x_{2} \\[1ex]
1 & 1 & 0 & 1 \\[1ex]
0 & 0 & 0 & 0
\end{array}
\right)
\]
 \[
A^{(4)}({\bf x})=
\left(
\begin {array}{cccc}
x_{2}x_{3} &   x_{1}x_{3} &  x_{1}x_{2} & 0\\[1ex]
-x_{2}-x_{3} &   -x_{1}-x_{3} &  -x_{1}-x_{2} & 0 \\[1ex]
1 & 1 & 1 & 0 \\[1ex]
0 & 0 & 0 & 0
\end{array}
\right).
\]
Es ist wiederum mit der Formel aus dem Satz von Vieta f\"ur $l\not=k$
\[
\Big( 1 \, \, t \, \, t^2 \cdots \, t^{n-1} \Big)\cdot
\Big( A^{(k)}({\bf x^{(m)}}) \Big)_{\cdot l} =  \displaystyle{
\prod \limits_{\begin{array}{c} \nu=1,\\ \nu \not=k, \\ \nu\not=l
\end{array}}^{n}  (t- x_{\nu}^{(m)}  ) }=:B^{(lk)}(t).
\]
Wir berechnen nun f\"ur $k=1,...,n$ die Matrizen
\[
F'({\bf x^{(m)}})^{-1} \cdot A^{(k)}({\bf x^{(m)}}).
\]
F\"ur das $(j,l)$-Element ergibt sich f\"ur $l\not= k$
\[
\frac{1}{B^{(j)}(x_{j}^{(m)})}
\Big( 1 \, \, x_{j}^{(m)} \, \, (x_{j}^{(m)})^2 \cdots \, (x_{j}^{(m)})^{n-1} \Big)\cdot
 \Big( A^{(k)}({\bf x^{(m)}}) \Big)_{\cdot l} =
\frac{B^{(lk)}(x_{j}^{(m)})}{B^{(j)}(x_{j}^{(m)})}
\]
und f\"ur $l= k$
\[
\frac{1}{B^{(j)}(x_{j}^{(m)})}
\Big( 1 \, \, x_{j}^{(m)} \, \, (x_{j}^{(m)})^2 \cdots \, (x_{j}^{(m)})^{n-1} \Big)\cdot
 o=0.
\]
Weiter ist f\"ur $l\not= k$
\[
\frac{B^{(lk)}(t)}{B^{(j)}(t)}=\left\{
\begin{array}{cl}
\displaystyle{\frac{t-x_{j}^{(m)}}{-\left(t-x_{l}^{(m)}\right)\cdot \left(t-x_{k}^{(m)}\right)} }
& \mbox{falls } j\not=l\mbox{ und } j\not=k,\\[4ex]
 \displaystyle{\frac{-1}{t-x_{k}^{(m)}}}  & \mbox{falls } j=l,\\[4ex]
  \displaystyle{\frac{-1}{t-x_{l}^{(m)}}}  & \mbox{falls } j=k.
\end{array}
\right.
\]
Somit gilt f\"ur $l\not= k$
\[
\frac{B^{(lk)}(x_{j}^{(m)})}{B^{(j)}(x_{j}^{(m)})}= \left\{
\begin{array}{cl}
0 & \mbox{falls } j\not=l\mbox{ und } j\not=k,\\[2ex]
 \displaystyle{\frac{1}{x_{k}^{(m)}-x_{j}^{(m)}}}  & \mbox{falls } j=l \mbox{ aber }l \not=k,\\[2ex]
 \displaystyle{\frac{1}{x_{l}^{(m)}-x_{j}^{(m)}}}  & \mbox{falls } j=k \mbox{ aber }l \not=k.
\end{array}
\right.
\]
Somit folgt $F'({\bf x^{(m)}})^{-1} \cdot A^{(k)}({\bf x^{(m)}})=$
\[
\left(
\begin{array}{ccccccc}
\displaystyle{\frac{1}{x_{k}^{(m)}-x_{1}^{(m)}}} & 0 & \cdots  & 0 & \cdots &\cdots & 0\\[2ex]
0 &  \displaystyle{\frac{1}{x_{k}^{(m)}-x_{2}^{(m)}}} & 0 & 0   &\cdots & \cdots &\vdots \\[2ex]
 0  & \ddots &  \ddots  &  0   &  \vdots  & & 0\\[2ex]
\displaystyle{\frac{1}{x_{1}^{(m)}-x_{k}^{(m)}}} & \cdots &
\displaystyle{\frac{1}{x_{k-1}^{(m)}-x_{k}^{(m)}}} & 0 &
\displaystyle{\frac{1}{x_{k+1}^{(m)}-x_{k}^{(m)}}} & \cdots &
\displaystyle{\frac{1}{x_{n}^{(m)}-x_{k}^{(m)}}} \\[2ex]
 0 & \cdots &
\cdots & 0 &
\displaystyle{\frac{1}{x_{k}^{(m)}-x_{k+1}^{(m)}}} & 0  &
0 \\[2ex]
\vdots  & \ddots &  \ddots  &  0   &    & \ddots & 0\\[2ex]
0 & \cdots & \cdots & 0 & \cdots & 0 & \displaystyle{\frac{1}{x_{k}^{(m)}-x_{n}^{(m)}}}
\end{array}
\right).
\]
Mit
\[
F'({\bf x^{(m)}})^{-1}F({\bf x^{(m)}}) =
\left(
\begin{array}{c}
\displaystyle{\frac{-p(x_{1}^{(m)})}{B^{(1)}(x_{1}^{(m)})}} \\[3ex]
\displaystyle{\frac{-p(x_{2}^{(m)})}{B^{(2)}(x_{2}^{(m)})}} \\[3ex]
\vdots \\[3ex]
\displaystyle{\frac{-p(x_{n}^{(m)})}{B^{(n)}(x_{n}^{(m)})}}
\end{array}
\right)
\]
folgt dann
\[
F'({\bf x^{(m)}})^{-1} \cdot A^{(k)}({\bf x^{(m)}}) \cdot \Big(
F'({\bf x^{(m)}})^{-1}  F({\bf x^{(m)}})
\Big)=\left(
\begin{array}{c}
\displaystyle{\frac{-p(x_{1}^{(m)})}{(x_{k}^{(m)}-x_{1}^{(m)}) \cdot B^{(1)}(x_{1}^{(m)})}} \\[3ex]
\vdots \\[3ex]
\displaystyle{\frac{-p(x_{k-1}^{(m)})}{(x_{k}^{(m)}-x_{k-1}^{(m)})\cdot B^{(k-1)}(x_{k-1}^{(m)})}} \\[5ex]
\displaystyle{\sum \limits_{\nu = 1, \, \nu \not=k}^{n}
\frac{-p(x_{\nu}^{(m)})}{(x_{\nu}^{(m)}-x_{k}^{(m)})\cdot B^{(\nu)}(x_{\nu}^{(m)})}} \\[5ex]
\displaystyle{\frac{-p(x_{k+1}^{(m)})}{(x_{k}^{(m)}-x_{k+1}^{(m)})\cdot B^{(k+1)}(x_{k+1}^{(m)})}} \\[3ex]
\vdots \\[3ex]
\displaystyle{\frac{-p(x_{n}^{(m)})}{(x_{k}^{(m)}-x_{n}^{(m)})
\cdot B^{(n)}(x_{n}^{(m)})}}
\end{array}
\right).
\]
Zuletzt folgt $F'({\bf x^{(m)}})^{-1} \cdot F''({\bf x^{(m)}}) \Big(
F'({\bf x^{(m)}})^{-1}  F({\bf x^{(m)}})\, , \, F'({\bf x^{(m)}})^{-1}  F({\bf x^{(m)}})
\Big)=$
\[
\left(
\begin{array}{c}
\displaystyle{\sum \limits_{ \nu = 1 \, \nu \not=1}^{n}
\frac{p(x_{\nu}^{(m)})\cdot p(x_{1}^{(m)})}{(x_{\nu}^{(m)}-x_{1}^{(m)})
B^{(\nu)}(x_{\nu}^{(m)})  B^{(1)}(x_{1}^{(m)})}} +
 \displaystyle{\sum \limits_{ k= 1 \, k \not=1}^{n}
\frac{p(x_{1}^{(m)})\cdot p(x_{k}^{(m)})}{(x_{k}^{(m)}-x_{1}^{(m)})
B^{(1)}(x_{1}^{(m)})  B^{(k)}(x_{k}^{(m)})}}\\[5ex]
\displaystyle{\sum \limits_{\nu = 1 \, \nu \not=2 }^{n}
\frac{p(x_{\nu}^{(m)})\cdot p(x_{2}^{(m)})}{(x_{\nu}^{(m)}-x_{2}^{(m)})
B^{(\nu)}(x_{\nu}^{(m)})  B^{(2)}(x_{2}^{(m)})}} +
 \displaystyle{\sum \limits_{k= 1 \, k \not=2 }^{n}
\frac{p(x_{2}^{(m)})\cdot p(x_{k}^{(m)})}{(x_{k}^{(m)}-x_{2}^{(m)})
B^{(2)}(x_{2}^{(m)})  B^{(k)}(x_{k}^{(m)})}}\\[5ex]
\vdots \\[5ex]
\displaystyle{\sum \limits_{\nu = 1 \, \nu \not=n }^{n}
\frac{p(x_{\nu}^{(m)})\cdot p(x_{n}^{(m)})}{(x_{\nu}^{(m)}-x_{n}^{(m)})
B^{(\nu)}(x_{\nu}^{(m)})  B^{(n)}(x_{n}^{(m)})}} +
 \displaystyle{\sum \limits_{  k= 1 \, k \not=n}^{n}
\frac{p(x_{n}^{(m)})\cdot p(x_{k}^{(m)})}{(x_{k}^{(m)}-x_{n}^{(m)})
B^{(n)}(x_{n}^{(m)})  B^{(k)}(x_{k}^{(m)})}}
\end{array}
\right)
\]
 \[
=2 \cdot \left(
\begin{array}{c}
\displaystyle{ \frac{p(x_{1}^{(m)})}{B^{(1)}(x_{1}^{(m)})} \sum \limits_{\nu = 1, \, \nu \not=1}^{n}
\frac{p(x_{\nu}^{(m)})}{(x_{\nu}^{(m)}-x_{1}^{(m)})\cdot
B^{(\nu)}(x_{\nu}^{(m)}) }}\\[5ex]
\displaystyle{ \frac{p(x_{2}^{(m)})}{B^{(2)}(x_{2}^{(m)})} \sum \limits_{\nu = 1, \, \nu \not=2}^{n}
\frac{p(x_{\nu}^{(m)})}{(x_{\nu}^{(m)}-x_{2}^{(m)})\cdot
B^{(\nu)}(x_{\nu}^{(m)}) }}\\[5ex]
\vdots \\[5ex]
\displaystyle{ \frac{p(x_{n}^{(m)})}{B^{(n)}(x_{n}^{(m)})} \sum \limits_{\nu = 1, \, \nu \not=n}^{n}
\frac{p(x_{\nu}^{(m)})}{(x_{\nu}^{(m)}-x_{n}^{(m)})\cdot
B^{(\nu)}(x_{\nu}^{(m)}) }}
\end{array}
\right).
\]
Man erh\"alt also das Tschebyschew-Verfahren
\[
{\bf x^{(m+1)}} =  {\bf x^{(m)}} -
\left(
\begin{array}{c}
\displaystyle{ \frac{p(x_{1}^{(m)})}{B^{(1)}(x_{1}^{(m)})}}
\left( \displaystyle{ \sum \limits_{\nu = 1, \, \nu \not=1}^{n}
\frac{p(x_{\nu}^{(m)})}{(x_{\nu}^{(m)}-x_{1}^{(m)})\cdot
B^{(\nu)}(x_{\nu}^{(m)})}-1} \right)\\[5ex]
\displaystyle{ \frac{p(x_{2}^{(m)})}{B^{(2)}(x_{2}^{(m)})}}
\left( \displaystyle{ \sum \limits_{\nu = 1, \, \nu \not=2}^{n}
\frac{p(x_{\nu}^{(m)})}{(x_{\nu}^{(m)}-x_{2}^{(m)})\cdot
B^{(\nu)}(x_{\nu}^{(m)})}-1 } \right)\\[5ex]
\vdots \\[5ex]
 \displaystyle{ \frac{p(x_{n}^{(m)})}{B^{(n)}(x_{n}^{(m)})}}
\left( \displaystyle{ \sum \limits_{\nu = 1, \, \nu \not=n}^{n}
\frac{p(x_{\nu}^{(m)})}{(x_{\nu}^{(m)}-x_{n}^{(m)})\cdot
B^{(\nu)}(x_{\nu}^{(m)})}-1 }\right)
\end{array}
\right).
\]
Mit (\ref{wfunk})  erh\"alt man dieselben Formeln wie in \cite{tanabe}. Siehe auch Formel (3.2) in \cite{kyurk}.
\section{Numerische Beispiele}
Bekanntlich ist das Newton-Verfahren lokal quadratisch konvergent, w\"ahrend die lokale Konvergenzordnung des Tschebyschew-Verfahrens 3 ist. Bei der praktischen Umsetzung bedeutet dies, dass sich bei quadratischer Konvergenz die Anzahl der exakten Stellen in jedem Iterationsschritt in etwa verdoppelt, w\"ahrend sie sich bei der Konvergenz\-ordnung 3 in etwa verdreifacht. F\"ur Fragen bzgl. globaler Konvergenz verweisen wir auf \cite{kyurk} und \cite{reinke}.
\begin{example} Wir betrachten
\[
p(t) = t^4 -5t^2 +6.
\]
$p(t)$ hat die Nullstellen
\[
\pm \sqrt{2}, \quad \pm \sqrt{3}.
\]
Wir w\"ahlen
\[
{\bf x^{(0)}}:=
\left(
\begin{array}{r} 1.2 \\1.8 \\ -1.2 \\-1.8 \end{array}
\right).
\]
Die ersten 5 Iterierten des Newton-Verfahrens f\"ur die positiven Wurzeln lauten
\begin{verbatim}
 1.402222222222222E+000   1.754074074074074E+000
 1.413432290193275E+000   1.732854607981912E+000
 1.414211612595975E+000   1.732052760484365E+000
 1.414213562361249E+000   1.732050807580748E+000
 1.414213562373095E+000   1.732050807568877E+000
\end{verbatim}
 Die ersten 5 Iterierten des Tschebyschew-Verfahrens f\"ur die positiven Wurzeln lauten
\begin{verbatim}
 1.403757613168724E+000   1.741105197378448E+000
 1.414197958229019E+000   1.732066406534148E+000
 1.414213562373021E+000   1.732050807568952E+000
 1.414213562373095E+000   1.732050807568877E+000
 1.414213562373095E+000   1.732050807568877E+000
\end{verbatim}
F\"ur die Wahl
\[
{\bf x^{(0)}}:=
\left(
\begin{array}{r} 1+i \\20+30i \\ 30+50i \\-40+30i \end{array}
\right)
\]
erhalten wir mit dem Abbruchkriterium $\| {\bf x^{(m+1)}}-{\bf x^{(m)}}\|_1 < 10^{-15}$ beim Newton-Verfahren
\[
\begin{array}{l}
 x_1^{(20)}=1.732050807568877 -5.392603844284260 \cdot 10^{-33} \cdot i\\
x_2^{(20)}= -1.414213562373095+  4.683861624749758 \cdot 10^{-31} \cdot i\\
x_3^{(20)}= 1.414213562373095+ 5.392603844284260 \cdot 10^{-33} \cdot i \\
x_4^{(20)}= -1.732050807568877 -7.395570986446986 \cdot 10^{-31} \cdot i
\end{array}
\]
und beim Tschebyschew-Verfahren
\[
\begin{array}{l}
x_1^{(16)}= 1.414213562373095+  0.000000000000000 \cdot 10^{0} \cdot i\\
x_2^{(16)}= -1.414213562373095 -8.407790785948902 \cdot 10^{-45} \cdot i\\
x_3^{(16)}= 1.732050807568877+  5.605193857299268 \cdot 10^{-45}  \cdot i \\
x_4^{(16)}=  -1.732050807568877+  4.203895392974451 \cdot 10^{-45} \cdot i
\end{array}
\]
\end{example}
\begin{example}
Um einen Eindruck zu erhalten, wie sich die Verfahren bei mehrfachen Nullstellen verhalten, betrachten 
wir  $p(t) = (t+1)^5= t^5 + 5 t^4 +10 t^3+10 t^2 + 5 t +1$.  Mit
\[
{\bf x^{(0)}}:=
\left(
\begin{array}{r} 1 \\2 \\ 3 \\4\\5 \end{array}
\right)
\]
und dem Abbruchkriterium $\| {\bf x^{(m+1)}}-{\bf x^{(m)}}\|_1 < 10^{-15}$ erhalten wir beim Newton-Verfahren
\[
\begin{array}{l}
 x_1^{(208)}= -0.9994349370863596  \\
x_2^{(208)}=   -1.000169058535995\\
x_3^{(208)}=  -0.9995475920797345 \\
x_4^{(208)}=  -1.000420616083746 \\
x_5^{(208)} = -1.000425875014050
\end{array}
\]
Beim Tschebyschew-Verfahren erh\"alt man
\[
\begin{array}{l}
x_1^{(316)}=  -1.000560780105459 \\
x_2^{(316)}=  -0. 9995856183790167  \\
x_3^{(316)}=  -1.000764338479257 \\
x_4^{(316)}=   -1.000765554071457 \\
x_5^{(316)} = -0.9993655607005603
\end{array}
\]
\end{example}


\begin{thebibliography}{99}
\bibitem{alefeld} G. Alefeld, J. Herzberger, {\it On the convergence speed of some algorithms for the simultaneous approximation of polynomial roots}, SIAM J. Numer. Anal., {\bf 11}  (1974), 237-243.
\bibitem{artin} E. Artin, {\it Galoissche Theorie}, Verlag Harri Deutsch, 1988.
\bibitem{bew} J. Bewersdorff, {\it Algebra f\"ur Einsteiger}, Vieweg, 2002.
\bibitem{froeba} S. Fr\"oba, A. Wassermann, {\it Die bedeutendsten Mathematiker}, Marix Verlag, Wiesbaden, 2007.
\bibitem{kanno} S. Kanno, N. V. Kjurkchiev, T. Yamamoto, {\it On some methods for the simultaneous determination of polynomial 
zeros}, Japan J. Indust. Appl. Math., {\bf 13} (1996), 267-288.
\bibitem{kerner} I. Kerner, {\it Ein Gesamtschrittverfahren zur
Berechnung der Nullstellen von Polynomen}, Numer. Math., \textbf{8} (1966), 290-294.
\bibitem{kyurk} N. V. Kyurkchiev, {\it Initial Approximations and Root Finding Methods}, Wiley-VCH Verlag Berlin, 1998. 
\bibitem{petkovic} M. S. Petkovic, {\it Iterative Methods for Simultaneous Inclusion of Polynomial Zeros}, Lecture Notes in Mathematics 1387, Springer-Verlag, Berlin, 1989.
\bibitem{reinke} B. Reinke, D. Schleicher, M. Stoll, {\it The Weierstra\ss-Durand-Kerner root finder is not generally convergent}, Math. Comp., {\bf 92} (2023), 839-866.
\bibitem{tanabe} K. Tanabe, {\it Behavior of the sequences around multiple zeros generated by some simultaneous methods for solving algebraic equations (Japanese)}, Tech. Rep. Inf. Process. Numer. Anal., 4-2 (1983), 1-6. 
\bibitem{weierstrass} K. Weierstra\ss, {\it Neuer Beweis des Satzes, dass jede ganze rationale Function einer Ver\"anderlichen dargestellt werden kann als ein Product aus linearen Functionen derselben Ver\"anderlichen}, in: Sitzungsberichte der k\"oniglich preussischen Akademie der Wissenschaften zu Berlin, 1891.
\bibitem{wloka} J. Wloka, {\it Funktionalanalysis und Anwendungen}, de Gruyter
Lehrbuch, 1971.
\end{thebibliography}
\end{document}